\documentclass[11pt,a4paper]{amsart}

\usepackage{amsmath,amsthm}

\theoremstyle{plain}
\newtheorem{thm}{Theorem}
\newtheorem{cor}[thm]{Corollary}

\usepackage{amsfonts}
\newcommand{\field}[1]{\mathbb{#1}}

\newcommand{\C}{\field{C}}

\newcommand{\loglike}[1]{\mathop{\mathrm {#1}}\nolimits}
\newcommand{\re}{\loglike{re}}

\newcommand{\diam}{\loglike{diam}}

\newcommand{\eps}{\varepsilon}
\newcommand{\calA}{{\mathcal A}}
\newcommand{\mycolon}{{:}\allowbreak\ }

\begin{document}

\title{Slices in the unit ball of a uniform algebra}

\author{Olav Nygaard and Dirk Werner}

\address{Department of Mathematics, Agder College, Tordenskjoldsgate 65, 
\qquad {} \linebreak 4604 Kristiansand, Norway}
\email{Olav.Nygaard@hia.no}

\address{Department of Mathematics, Freie Universit\"at Berlin,
Arnimallee~2--6, \qquad {}\linebreak D-14\,195~Berlin, Germany}
\email{werner@math.fu-berlin.de}


\date{\today}

\begin{abstract} We show that every nonvoid relatively weakly open subset,
in particular every slice, of the unit ball of an 
infinite-dimensional uniform algebra has diameter~$2$.
\end{abstract}

\maketitle

It is an important task in Banach space theory to determine the extreme 
point structure of the unit ball for various examples of Banach spaces. 
The most common way to describe ``corners" of convex sets is by 
looking for \em extreme points, exposed points, denting points\/ 
\em and \em strongly exposed points\em. Every strongly exposed point 
is both denting and exposed and every denting or exposed point is extreme. 

In this note $\calA$ denotes an infinite-dimensional uniform algebra,
i.e., an infinite-dimensional closed subalgebra of some $C(K)$-space
which separates the points of $K$ and contains the constant functions. 
In \cite{benwie} Beneker and Wiegerinck demonstrated the 
non-existence of strongly exposed points in $B_{\calA}$, 
the closed unit ball of $\calA$. 
Here we shall prove a stronger result by more elementary means.
A corollary of our result is that the set of denting points is, 
in fact, also empty.

Recall that we may assume that $K$ is the Silov boundary of $\calA$.
It is a fundamental result in the theory of uniform algebras that
then the set of strong boundary points is dense in $K$; cf.\
\cite[p.~48 and p.~78]{leibo}. 
(A point $x\in K$ is a strong boundary point if for every
neighbourhood $V$ of $x$ and every $\delta>0$
there is some $f\in\calA$ such that
$f(x)=\|f\|=1$ and $|f|\le\delta$ off $V$.)

We now turn to our  first result, which gives a quantitative statement 
of the non-dentability of $B_\calA$.

\begin{thm}\label{thm1}
Every slice of the unit ball of an infinite-dimensional 
uniform algebra $\calA$ has diameter~$2$. 
\end{thm}

\begin{proof} 
Take an arbitrary slice
$S=\{a\in B_{\calA}\mycolon \re \ell(a)\ge1-\eps\}$, where $\|\ell\|=1$. 
We will produce two functions in $S$ having distance nearly~$2$.

Let $0<\delta\le\eps/11$. We first pick some $f\in B_{\calA}$ such that
$$
\re \ell(f)\ge 1-\delta.
$$
The functional $\ell$ can be represented by a regular Borel measure
$\mu$ on $K$ with $\|\mu\|=1$,
i.e.,  $\ell(a)=\int_K a\,d\mu$ for all $a\in\calA$. 
Let $\emptyset \neq V_{0}\subset K$ be
an open set with $|\mu|(V_{0})\le\delta$; such a set exists since $K$ is
infinite. Fix a strong boundary point $x_{0}\in V_{0}$. Using 
the definition of a strong boundary point, inductively construct functions
$g_{1},g_{2},\ldots\in \calA$ and nonvoid open subsets $V_{0}\supset
V_{1} \supset V_{2} \supset \ldots$      such that
$$
g_{n}(x_{0})= \|g_{n}\| =1,\ \ |g_{n}|\le\delta \mbox{ on }K\setminus V_{n-1}
$$
and 
$$
V_{n}=\{x\in V_{n-1}\mycolon |g_{n}(x)-1| <\delta\}.
$$
Let $N>1/\delta$ and define
$$
g= \frac1N \sum_{k=1}^N g_{k},\ \ h=f(1-g)\in \calA.
$$

By construction, $|h|\le\delta$ on $V_{N}$ and $|h|\le1+\delta$ on
$K\setminus V_{0}$. We claim that $\|h\|\le 1+3\delta$. In fact, if
$x\in V_{r-1}\setminus V_{r}$, then 
$|1-g_{k}(x)|\le\delta$ if $1\le k<r$,
$|g_{r}(x)|\le1$ and
$|g_{k}(x)|\le\delta$ if $r<k\le N$, and therefore
$$
|h(x)|\le
\frac1N \sum_{k=1}^N |1-g_{k}(x)| \le \frac{(N-1)(1+\delta)+2}N
\le 1+3\delta.
$$

We now estimate $|\ell(f)-\ell(h)|$:
\begin{align*}
|\ell(f)-\ell(h)| &\le
\int_{K\setminus V_{0}} |f-h|\,d|\mu| +
\int_{ V_{0}}   |f-h|\,d|\mu| \\
&\le
\int_{K\setminus V_{0}} |g|\,d|\mu| +
\int_{ V_{0}}   (|f|+|h|)\,d|\mu| \\
&\le
\delta + (2+3\delta)|\mu|(V_{0}) \le 4\delta.
\end{align*}

Next, we produce a function $\varphi\in \calA$
such that
$$
\varphi(x_{0})= \|\varphi\| =1,\ \ |\varphi|\le\delta \mbox{ on }
K\setminus V_{N}.
$$
We then have $\|h\pm\varphi\|\le 1+4\delta$, and the functions 
$\psi_{\pm} = (h\pm\varphi)/(1+4\delta)$ are in the unit ball of $\calA$.
We have
$
|\ell(\varphi)|\le
 2\delta
$
and thus
$$
|\ell(\psi_{\pm}) - \ell(h)|  \le
|\ell(h)| \frac{4\delta}{1+4\delta} + \frac{2\delta}{1+4\delta}
\le 6\delta.
$$
Consequently,
$$
\re \ell(\psi_{\pm}) \ge \re \ell(f) -10\delta \ge 1-11\delta \ge
1-\eps
$$
so that $\psi_{\pm}\in S$; but
$\|\psi_{+}-\psi_{-}\|=2/(1+4\delta)\to2$ as $\delta\to0$. Hence
$\diam S=2$.
\end{proof}

The point of working with $g$  rather than $g_{1}$ in the 
proof is to control $\|1-g\|$. Another way to achieve this is to
construct a suitable conformal map $\phi$ from the unit disk to a
neighbourhood of $[0,1]$ in $\C$ and to consider $\phi\circ g_{1}$.

We now extend Theorem~\ref{thm1} to relatively weakly open subsets.

\begin{thm}\label{thm2}
Every nonvoid relatively weakly open subset $W$
of the unit ball of an infinite-dimensional 
uniform algebra $\calA$ has diameter~$2$. 
\end{thm}

\begin{proof}
Every nonvoid relatively weakly open subset 
of the unit ball of a Banach space contains a convex combination of slices, 
see \cite[Lemma~II.1]{GGMS} or \cite{Shv}. Thus, if $W\subset B_{\calA}$
is given as above, there are slices $S^{(1)},\dots,S^{(n)}$ and $0\le 
\lambda_{j}\le1$, $\sum_{j=1}^n \lambda_{j}=1$, such that $\sum_{j=1}^n
\lambda_{j}S^{(j)}\subset W$.

Let $S^{(j)}=\{a\in B_{\calA}\mycolon 
\re \ell_{j}(a) \ge 1-\eps_{j}\}$ with $\|\ell_{j}\|
=1$ and representing measures $\mu_{j}$. We now perform the construction 
of the proof of Theorem~\ref{thm1} with $\eps= \min\eps_{j}$, $0<\delta\le
\eps/11$ as before and a nonvoid open set $V_{0}\subset K$ such that
$|\mu_{j}|(V_{0})\le\delta$ for all $j$. We obtain functions $h^{(j)}$ and 
$\varphi$ (independently of~$j$) such that $(h^{(j)}\pm \varphi)/(1+4\delta)
\in S^{(j)}$ and $\|\varphi\|=1$. Therefore $\sum_{j=1}^n \lambda_{j}h^{(j)}
\pm \varphi \in (1+4\delta)W$, and $\diam W=2$.
\end{proof}

In  case $K$ does not have isolated points, Theorem~\ref{thm1} is a
formal consequence of the Daugavet property of $\calA$, proved in 
\cite{Woj92} or \cite{Dirk10}, and \cite[Lemma~2.1]{KSSW}. Likewise 
Theorem~\ref{thm2} follows from \cite{Shv}.

\begin{cor}\label{cor3}
The unit ball of an infinite-dimensional uniform algebra does not contain
any denting points or merely points of continuity for the identity mapping
with respect to the weak and the norm topology.
\end{cor}

We would also like to remark that
T.S.S.R.K. Rao \cite{rao} has shown that $B_{WC(K,X)}$, 
the unit ball in the space of continuous functions from a 
compact Hausdorff space into a Banach space equipped with its 
weak topology, has no denting points.
He has also given a proof of Corollary~\ref{cor3} based on techniques from 
that paper.


\end{document}